Some Results on the Counterfeit Coins Problem

Li An-Ping

Beijing 100085, P.R.China
apli0001@sina.com

Abstract

We will present some results on the counterfeit coins problem in the case of multi-sets.

## I. Introduction

Searching the counterfeit coins from the given coins by a balance is a well-known problem of combinatorial search, which there are several versions and a intensive researches, for the detail to see the papers [1]~[9].

From this paper, we will discuss some the called multi-sets cases. Let $S_1, \cdots, S_m$ be $m$ sets of coins, in which each set $S_i$ just contains one false coin, and the fakes are known with same weight but heavier than the normals, where $S_i \cap S_j = \varnothing$, for all $i \neq j$. Suppose that $|S_i| = n_i$, $1 \leq i \leq m$, define $g_1(n_1, n_2, \cdots, n_m)$ the least number of weighings to find all the $m$ fakes in $\bigcup_{1 \leq i \leq m} S_i$, which will be simply written as $g_1(n|_m)$ when $n_i = n, 1 \leq i \leq m$. The present paper will give some estimations for $g_1(n|_m)$.

In the rest of the section, we introduce some symbols and notations which will be used in this paper.

$L:R$ : A weighing or a comparison of coin set $L$ against coin set $R$.

$L > R$ ( $L < R$, $L = R$ ) : $L$ is heavier ( lighter, equal to) than $R$ respectively.

$X \cdot Y$ : Cartesion product of sets $X$ and $Y$.

$\prec$ : A partial order of the space $\bigcup_{k=0}^{\infty} X^k$, $\forall \alpha, \beta \in \bigcup_{k=0}^{\infty} X^k$, $\alpha \prec \beta \Leftrightarrow \alpha = (x_1, \cdots, x_l)$, and $\beta = (x_1, \cdots, x_l, \cdots)$.

$X_{\{i,j,\cdots,l\}} := X_i \cup X_j \cup \cdots \cup X_l$, and simply written as $X_{i,j,\cdots,l}$ if no confusability.

$\lceil x \rceil$ : The least integer no less than the real number $x$.

$Ct(X)$ : The number of fakes in the set $X$.

Let $A_0 = \{0\}$, $A = \{+, 0, -\}$, for a positive integer $k$, we define $A_k = A_0 \cdot A^{k-1}$, $\mathcal{A} = \bigcup_{k=1}^{\infty} A_k$, and $\mathcal{W} = \{ L:R \mid L, R \subset \bigcup_{i=1}^{m} S_i, |L| = |R| \}$.

It is clear that an algorithm of searching for the counterfeit coins is just a map $F : \mathcal{A} \to \mathcal{W}$.

Let $\mathcal{S} = \prod_{i=1}^{m} S_i$, then a algorithm $F$ will induces a map $f: \mathcal{S} \to \mathcal{A}$: $X \mapsto (0, a_1, a_2, \cdots)$,

where

$a_1 = \text{sgn}(|L_1 \cap X| - |R_1 \cap X|)$, $F(0) = L_1 : R_1$,
$a_i = \text{sgn}(|L_i \cap X| - |R_i \cap X|)$, $F((0, a_1, \cdots, a_{i-1})) = L_i : R_i$, for $i > 1$.

Hence, $\forall \alpha \in \mathcal{A}$, $\alpha$ determines a subset $\mathcal{S}_\alpha$ of $\mathcal{S}$: $\mathcal{S}_\alpha = \{P | P \in \mathcal{S}, \alpha \prec f(P)\}$, we call $\mathcal{S}_\alpha$ as the objective set on the direction $\alpha$.

**Definition 1.** For a algorithm $F$ and a positive integer $k$, if $\forall \alpha \in \mathcal{A}_k, |\mathcal{S}_\alpha| \leq 1$, then algorithm $F$ is called $k$-$completed$.

## 2. Main results

It is well-known that

$$g_1(n) = \lceil \log_3 n \rceil. \tag{2.1}$$

So, there is a simple estimation

$$\left\lceil \sum_{i=1}^{m} \log_3 n_i \right\rceil \leq g_1(n_1, \cdots, n_m) \leq \sum_{i=1}^{m} \lceil \log_3 n_i \rceil. \tag{2.2}$$

In this paper, our main results are following estimate

**Proposition 1**

$$g_1(n|_2) \leq \lceil \log_3 n \rceil + \lceil \log_3(n/5) \rceil + 1, \tag{2.3}$$

$$g_1(n|_3) \leq \lceil \log_3 n \rceil + \lceil \log_3(n/4) \rceil + \lceil \log_3(n/6) \rceil + 2, \tag{2.4}$$

$$g_1(n|_4) \leq \lceil \log_3 n \rceil + \lceil \log_3(n/5) \rceil + \lceil \log_3(n/11) \rceil + \lceil \log_3(n/20) \rceil + 5, \tag{2.5}$$

$$g_1(n|_5) \leq \lceil \log_3 n \rceil + \lceil \log_3(n/7) \rceil + \lceil \log_3(n/11) \rceil + \lceil \log_3(n/13) \rceil + \lceil \log_3(n/17) \rceil + 7. \tag{2.6}$$

To prove the result above need the following individual results

**Lemma 1**

$$g_1(2|_3) = 2, \tag{2.7}$$
$$g_1(2, 4) = 2, \tag{2.8}$$
$$g_1(4, 20) = 4, \tag{2.9}$$

$$g_1(5|_2) = 3, \tag{2.10}$$
$$g_1(11|_4) = 9, \tag{2.11}$$
$$g_1(7|_5) = 9, \tag{2.12}$$
$$g_1(11|_5) = 11, \tag{2.13}$$
$$g_1(13|_5) = 12, \tag{2.14}$$
$$g_1(17|_5) = 13. \tag{2.15}$$

Proof. The identities (2.7) and (2.8) are easier and left to the readers. The algorithms for the rest are put in the end as an appendix.

The results above have the following corollary

**Corollary 1**

$$g_1(4|_3) = 4, \tag{2.16}$$
$$g_1(20|_4) = 11. \tag{2.17}$$

Proof. (2.12) is followed by applying (2.7) twice, and (2.13) is from (2.9) and (2.10).

Proof of Proposition 1. Let

$$n = \lambda \cdot 3^k, \quad 1 \leq \lambda \leq 5,$$

Where $\lambda$ is a real number and $k$ is a non-negative integer. We will take induction on $k$. By (2.6) of Lemma 1, (2.3) is stand for $k=0$. For $k>0$, suppose that $A$ and $B$ are two coin sets, $|A|=|B|=n, Ct(A)=Ct(B)=1$. Take subsets $A_1, A_2 \subset A$, and $B_1, B_2 \subset B$, such that $|A_1|=|A_2|=|B_1|=|B_2|=\lceil n/3 \rceil$, and make two weighings $A_1 : A_2$ and $B_1 : B_2$, then it will be found the subsets $A' \subset A$ and $B' \subset B$ with $|A'| \leq \lceil n/3 \rceil$, $|B'| \leq \lceil n/3 \rceil$, $Ct(A')=1$, $Ct(B')=1$. Without loss generality, we may assume that $|A'|=|B'|$, then by the induction,

$$g_1(n|_2) \leq 2 + \lceil \log_3 \lceil n/3 \rceil \rceil + \lceil \log_3 \lceil n/5 \rceil \rceil + 1$$
$$\leq \lceil \log_3 n \rceil + \lceil \log_3(n/5) \rceil + 1.$$

The proofs for (2.4), (2.5) and (2.6) are similar but instead apply (2.7), (2.12) and (2.10), (2.11), (2.17) and (2.13), (2.14), (2.15), (2.16) respectively, so which are omitted.

*Note.* The upper bounds in Proposition 1 may be written as following form

$$g_1(n|_2) \leq \lceil 2 \cdot \log_3 n + 0.071 \rceil, \qquad (2.3')$$
$$g_1(n|_3) \leq \lceil 3 \cdot \log_3 n + 0.218 \rceil, \qquad (2.4')$$
$$g_1(n|_4) \leq \lceil 4 \cdot \log_3 n + 0.274 \rceil. \qquad (2.5')$$
$$g_1(n|_5) \leq \lceil 5 \cdot \log_3 n + 0.325 \rceil. \qquad (2.6')$$

**Appendix**    Algorithms

The sketch of algorithm $g_1(4, 20) = 4$

Let $X = \{x_i \mid 1 \le i \le 4\}$, $Y = \{y_i \mid 1 \le i \le 20\}$, $Ct(X) = Ct(Y) = 1$. The following is a $4\text{-}completed$ search algorithm

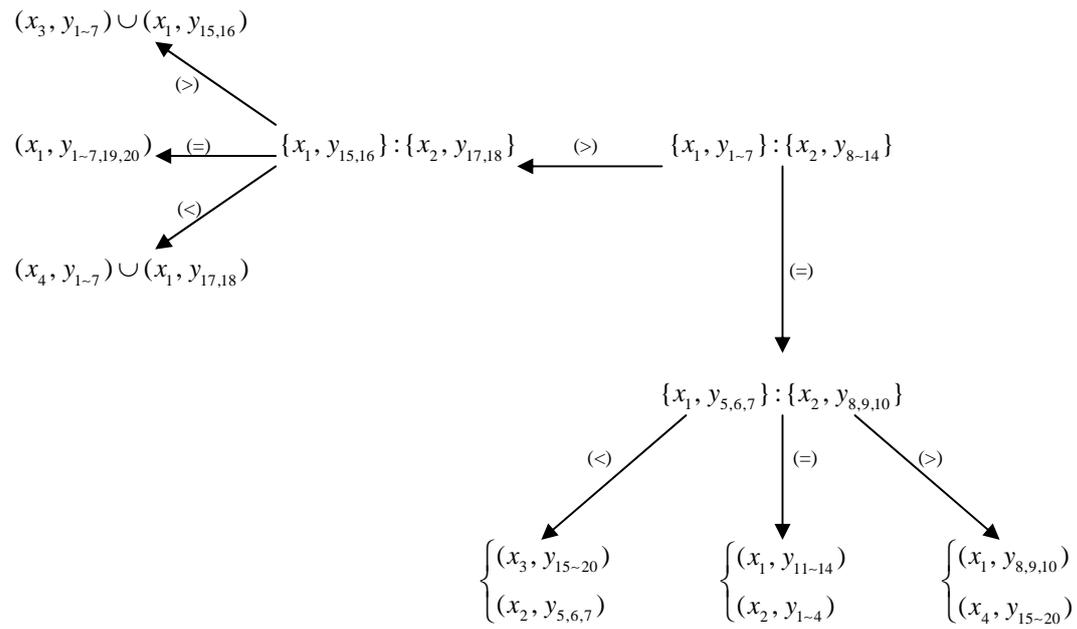

Fig.1

The rest two weighings are a routine work, so omitted.

The sketch of algorithm $g_1(5|_2) = 3$

Let $A = \{a_i | 1 \leq i \leq 5\}$, $B = \{b_i | 1 \leq i \leq 5\}$, $Ct(A) = Ct(B) = 1$. A $3$-*completed* algorithm is as following

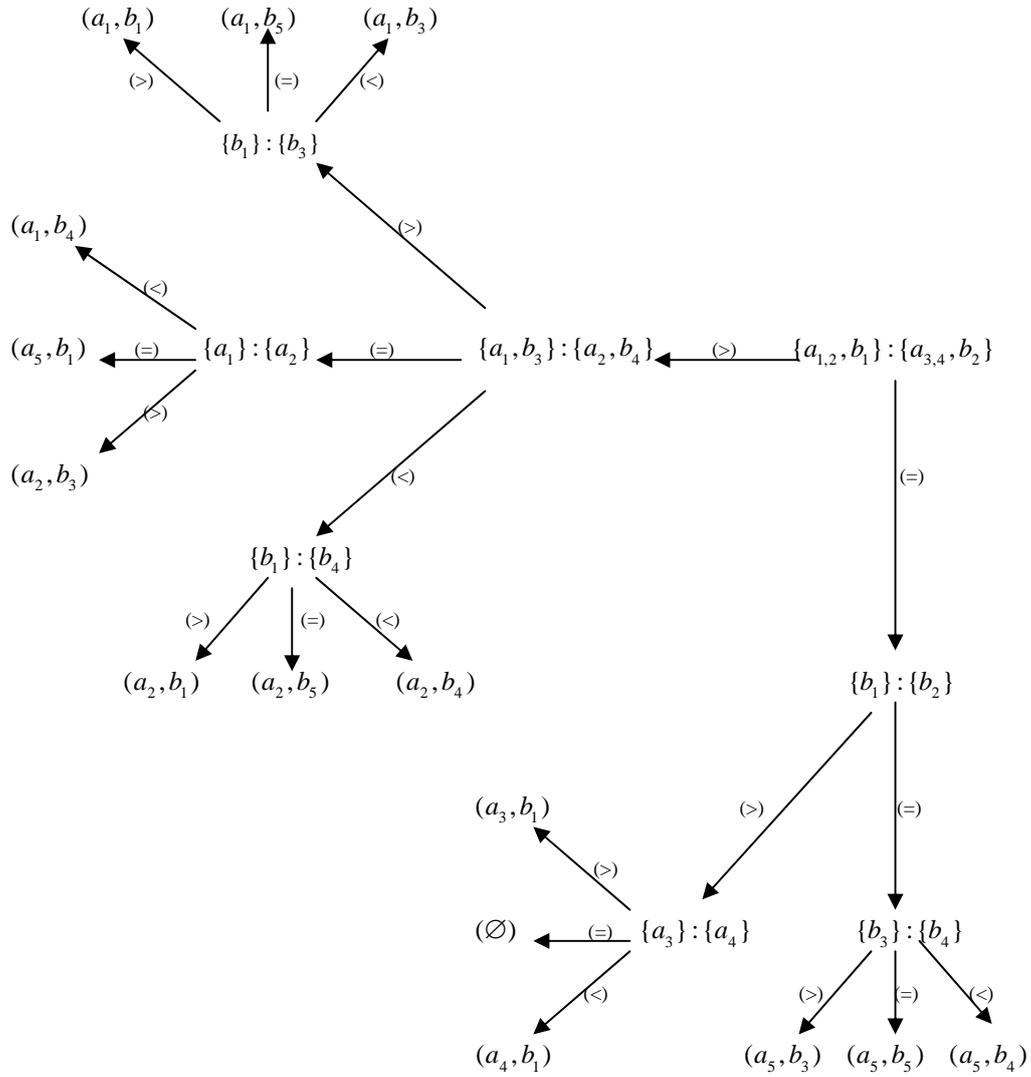

Fig.2

The sketch of algorithm $g_1(11|_4) = 9$

Let $A = \{a_i \mid 1 \le i \le 11\}$, $B = \{b_i \mid 1 \le i \le 11\}$, $C = \{c_i \mid 1 \le i \le 11\}$, $D = \{d_i \mid 1 \le i \le 11\}$, $Ct(A) = Ct(B) = Ct(C) = Ct(D) = 1$. For $X = A, B, C,$ or $D$, let $X = \bigcup_{0 \le i \le 3} X_i$, $|X_i| = 3$, $i = 1, 2, 3$, and $|X_0| = 2$. The sketch of a feasible algorithms is as following

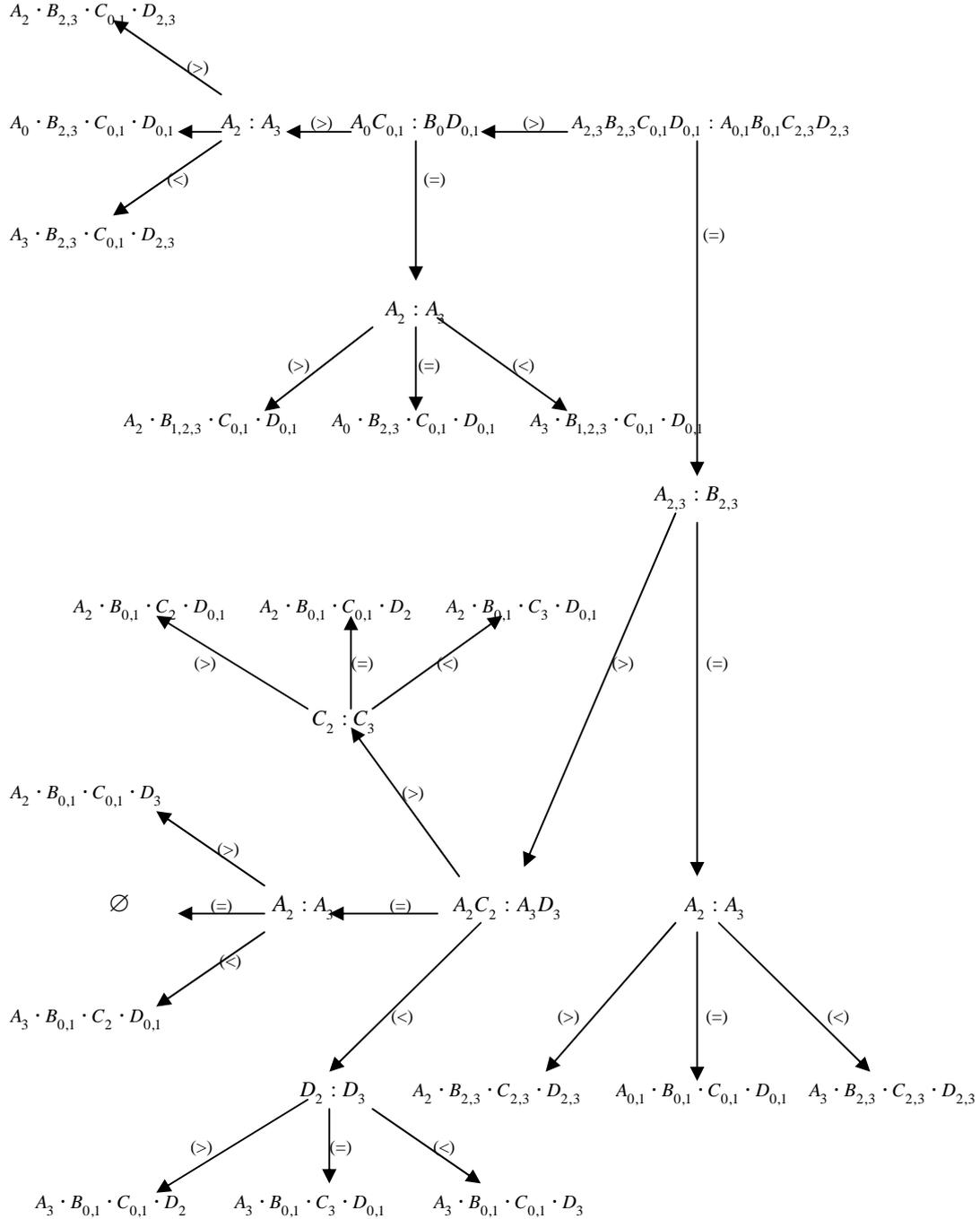

Fig.3

The sketch of algorithm $g_1(11|_5) = 11$

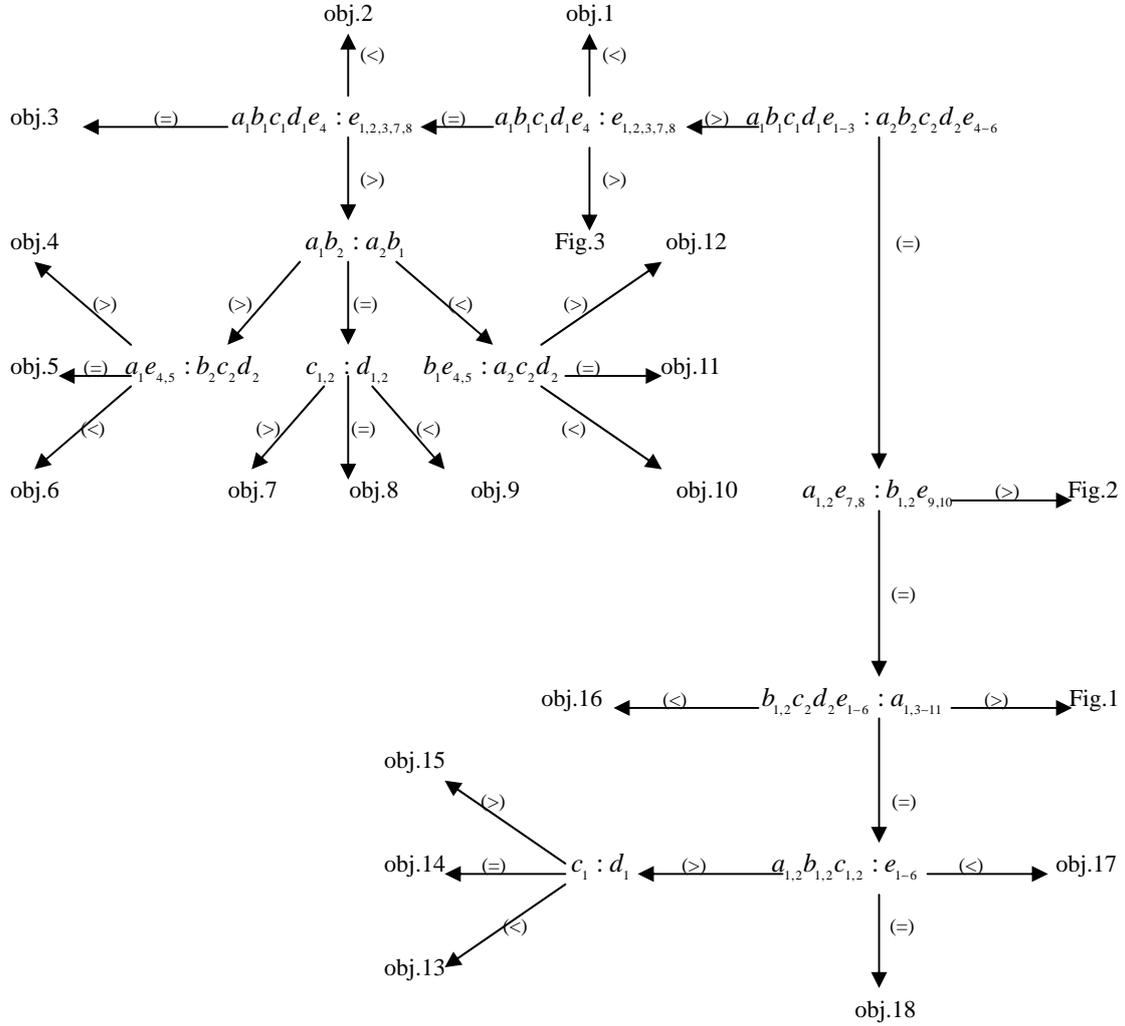

Fig.A

$$obj.1 = a_{3-11} \cdot b_{3-11} \cdot c_{3-11} \cdot d_{3-11} \cdot e_{1,2,3}, \qquad obj.4 = a_1 \cdot b_{3-11} \cdot c_{3-11} \cdot d_{3-11},$$

$$obj.2 = e_{7,8} \cdot \begin{cases} a_1 \cdot b_{3-11} \cdot c_{3-11} \cdot d_{3-11} \\ b_1 \cdot a_{3-11} \cdot c_{3-11} \cdot d_{3-11} \\ c_1 \cdot a_{3-11} \cdot b_{3-11} \cdot d_{3-11} \\ d_1 \cdot a_{3-11} \cdot b_{3-11} \cdot c_{3-11} \end{cases}, \qquad obj.3 = e_{1,3} \cdot \begin{cases} a_1 \cdot b_{3-11} \cdot c_{3-11} \cdot d_{3-11} \\ b_1 \cdot a_{3-11} \cdot c_{3-11} \cdot d_{3-11} \\ c_1 \cdot a_{3-11} \cdot b_{3-11} \cdot d_{3-11} \\ d_1 \cdot a_{3-11} \cdot b_{3-11} \cdot c_{3-11} \end{cases}$$

$$obj.5 = \begin{cases} a_1 \cdot b_2 \cdot e_{1-3} \cdot c_{3-11} \cdot d_{3-11} \\ a_1 \cdot c_2 \cdot e_{1-3} \cdot b_{3-11} \cdot d_{3-11} \\ a_1 \cdot d_2 \cdot e_{1-3} \cdot b_{3-11} \cdot d_{3-11} \end{cases}, \qquad obj.6 = \begin{cases} b_2 \cdot c_1 \cdot e_{1-3} \cdot a_{3-11} \cdot d_{3-11} \\ b_2 \cdot d_1 \cdot e_{1-3} \cdot a_{3-11} \cdot c_{3-11} \end{cases}$$

$$obj.7 = c_1 \cdot e_2 \cdot a_{3-11} \cdot b_{3-11} \cdot d_{3-11}, \qquad obj.9 = d_1 \cdot e_2 \cdot a_{3-11} \cdot b_{3-11} \cdot c_{3-11}$$

$$obj.8 = \begin{cases} c_1 \cdot d_2 \cdot e_{1-3} \cdot a_{3-11} \cdot b_{3-11} \\ c_2 \cdot d_2 \cdot e_{1-3} \cdot a_{3-11} \cdot b_{3-11} \end{cases}, \qquad obj.10 = \begin{cases} a_2 \cdot c_1 \cdot e_{1-3} \cdot b_{3-11} \cdot d_{3-11} \\ a_2 \cdot d_1 \cdot e_{1-3} \cdot b_{3-11} \cdot c_{3-11} \end{cases}$$

$$obj.11 = b_1 \cdot e_{1-3} \cdot \begin{cases} a_2 \cdot c_{3-11} \cdot d_{3-11} \\ c_2 \cdot a_{3-11} \cdot d_{3-11} \\ d_2 \cdot a_{3-11} \cdot c_{3-11} \end{cases}, \qquad obj.12 = b_1 \cdot e_2 \cdot a_{3-11} \cdot c_{3-11} \cdot d_{3-11},$$

$$obj.13 = \begin{cases} a_2 \cdot d_1 \cdot e_{9,10} \cdot b_{3-11} \cdot c_{3-11} \\ b_2 \cdot d_1 \cdot e_{7,8} \cdot a_{3-11} \cdot c_{3-11} \\ c_2 \cdot d_1 \cdot e_{11} \cdot a_{3-11} \cdot b_{3-11} \end{cases}, \qquad obj.14 = \begin{cases} a_1 \cdot c_2 \cdot e_{9,10} \cdot b_{3-11} \cdot d_{3-11} \\ a_1 \cdot d_2 \cdot e_{9,10} \cdot b_{3-11} \cdot c_{3-11} \\ a_1 \cdot b_2 \cdot e_{11} \cdot c_{3-11} \cdot d_{3-11} \end{cases}$$

$$obj.15 = \begin{cases} a_2 \cdot c_1 \cdot e_{9,10} \cdot b_{3-11} \cdot c_{3-11} \\ b_2 \cdot c_1 \cdot e_{7,8} \cdot a_{3-11} \cdot c_{3-11} \\ d_2 \cdot c_1 \cdot e_{11} \cdot a_{3-11} \cdot b_{3-11} \end{cases}, \qquad obj.16 = e_{11} \cdot a_{3-11} \cdot b_{3-11} \cdot c_{3-11} \cdot d_{3-11},$$

$$obj.17 = c_1 \cdot e_{4,5,6} \cdot a_{3-11} \cdot b_{3-11} \cdot d_{3-11}, \qquad obj.18 = d_1 \cdot e_{4,5,6} \cdot a_{3-11} \cdot b_{3-11} \cdot c_{3-11}$$

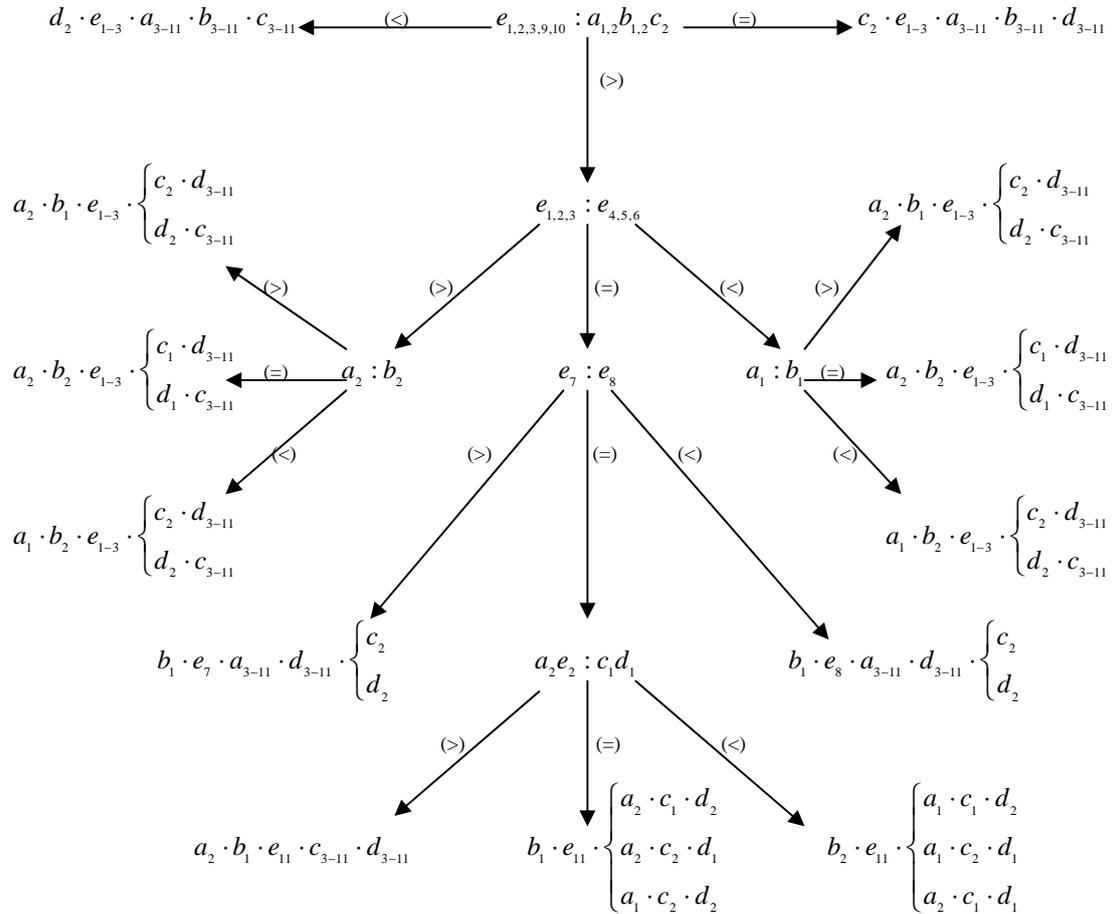

Fig.1

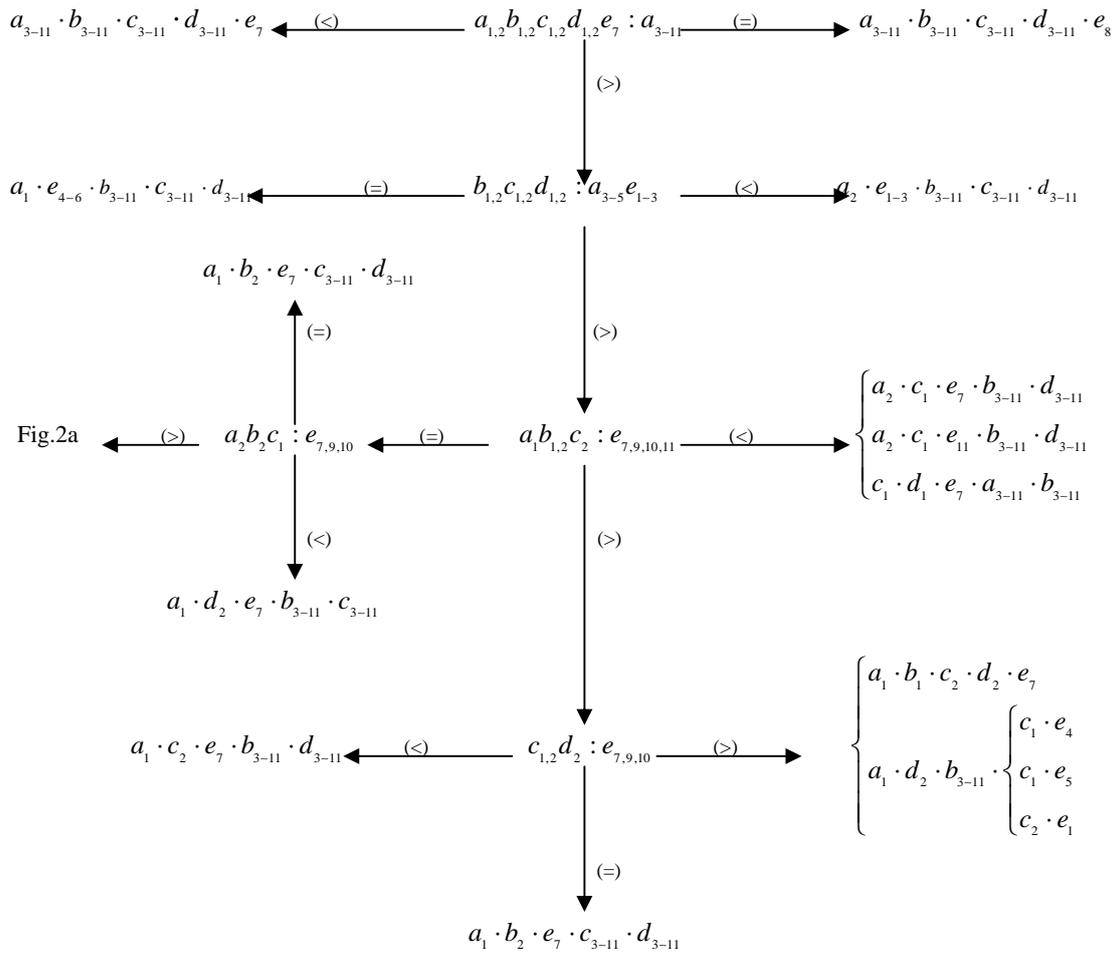

Fig. 2

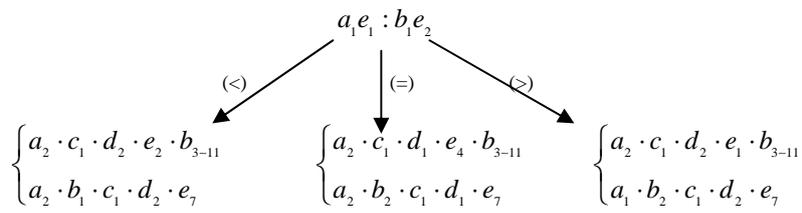

Fig.2a

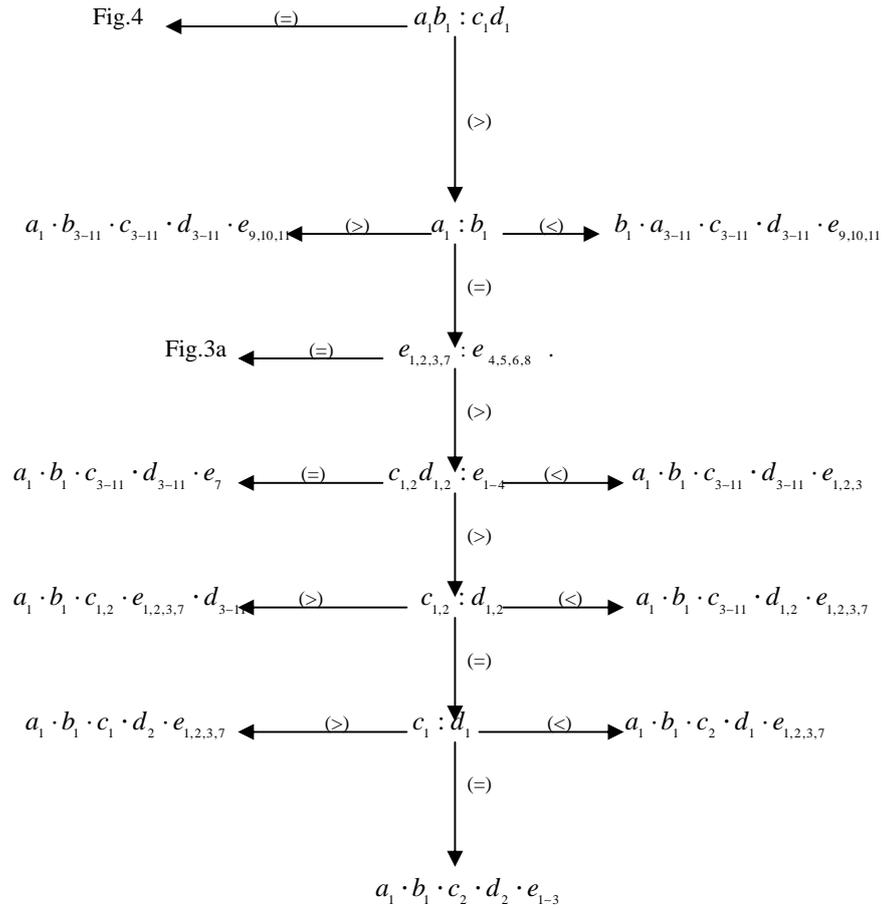

Fig.3

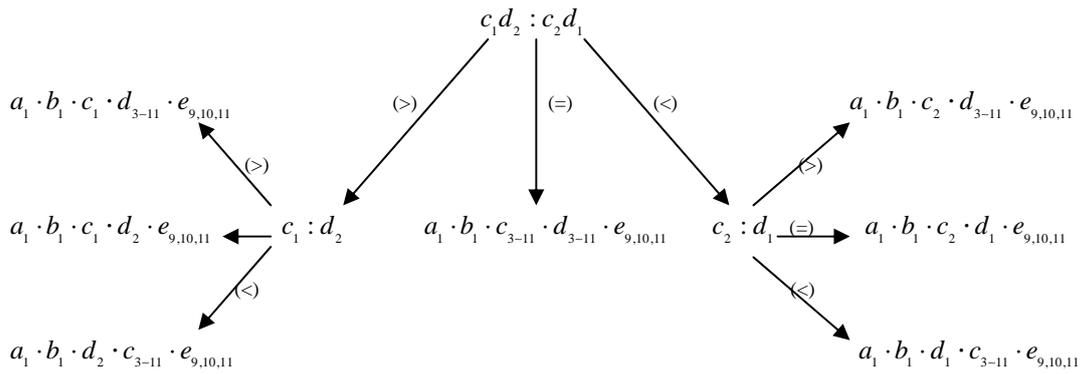

Fig. 3a

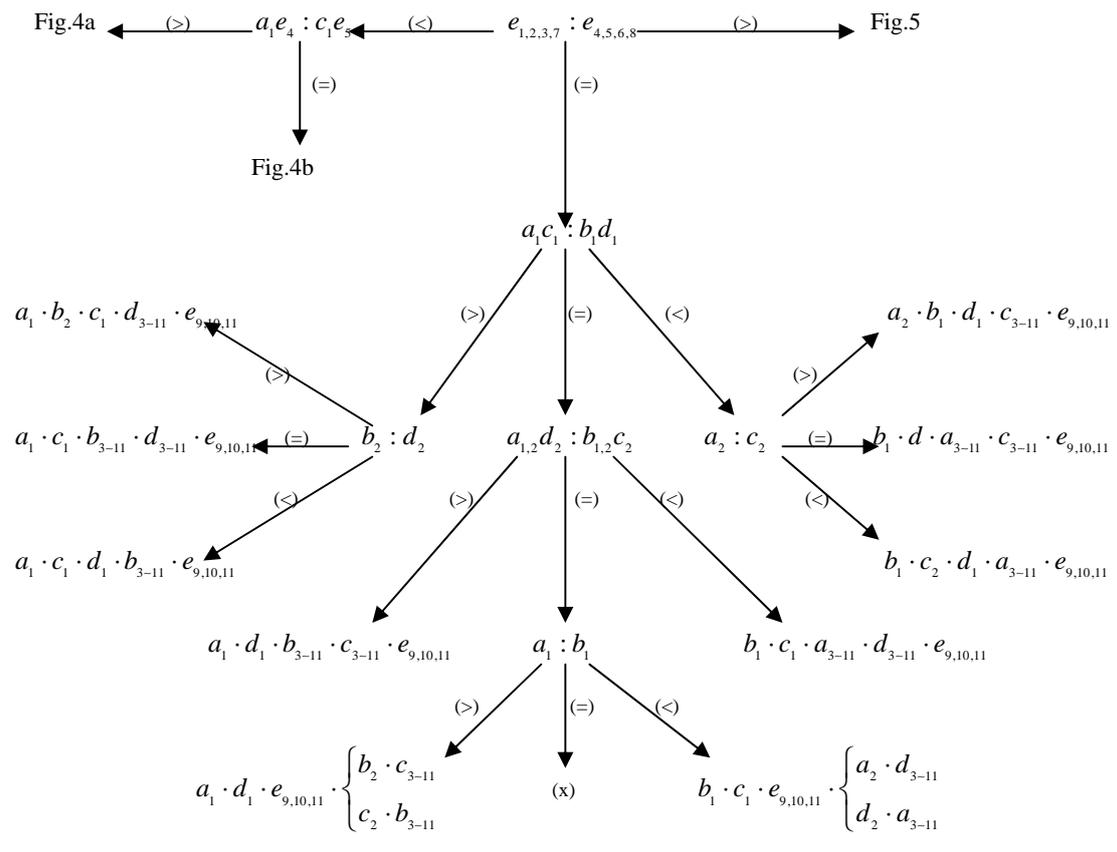

Fig.4

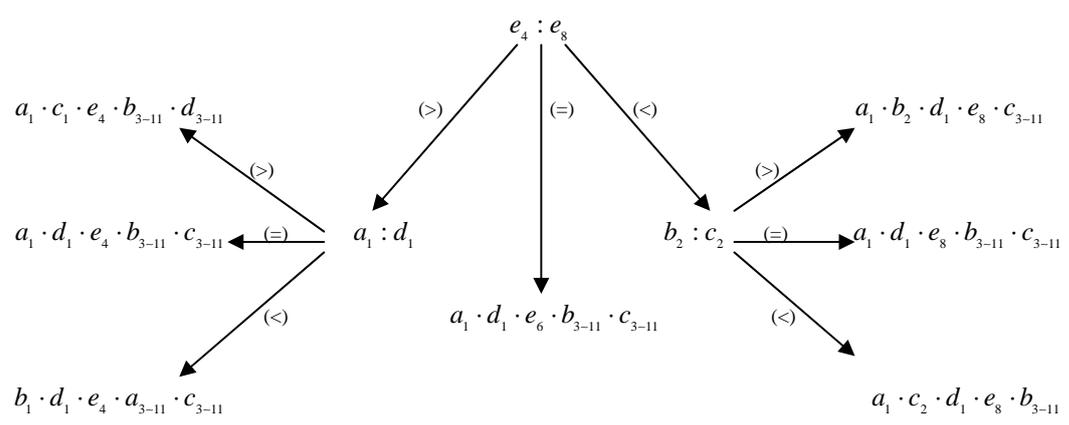

Fig.4a

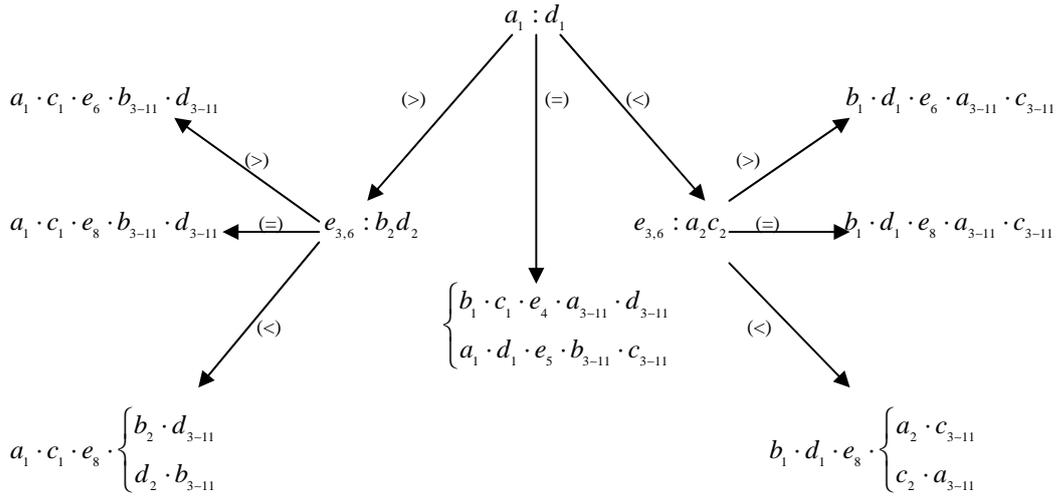

Fig.4b

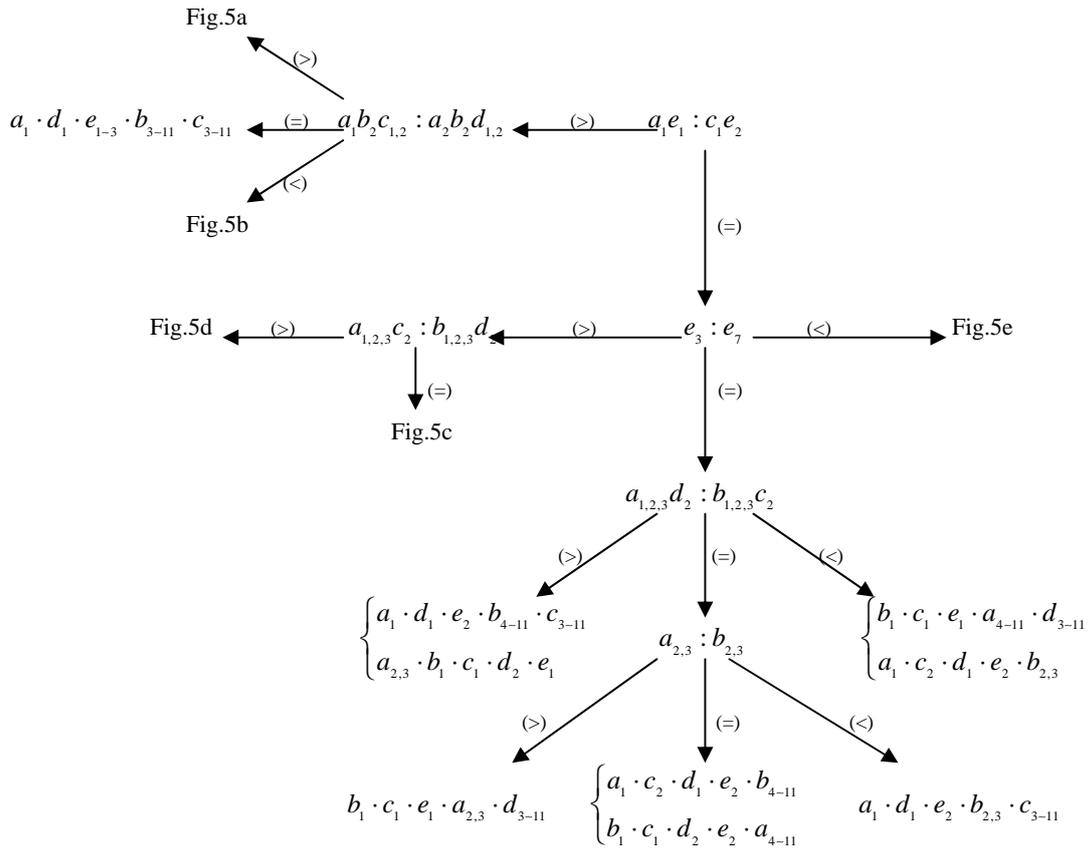

Fig.5

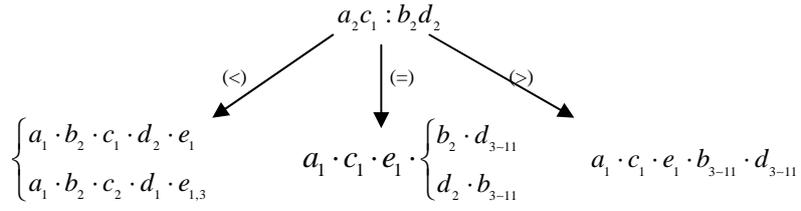

Fig.5a

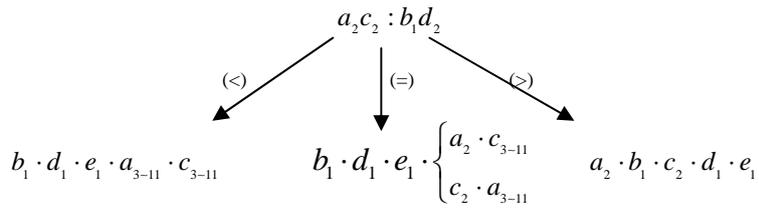

Fig.5b

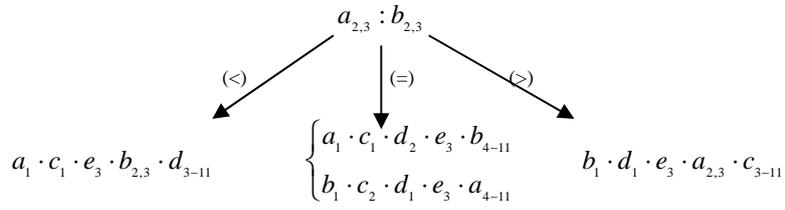

Fig.5c

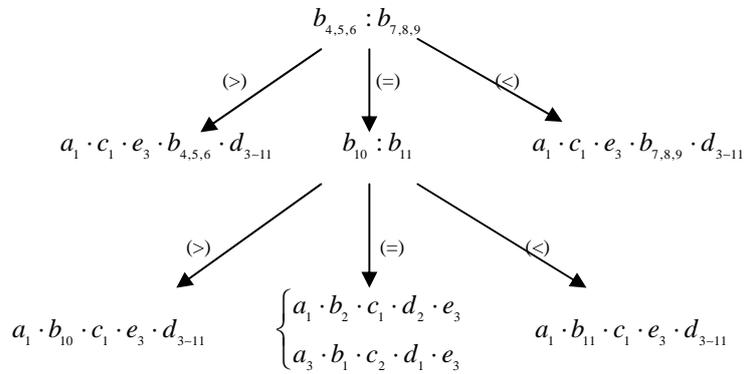

Fig.5d

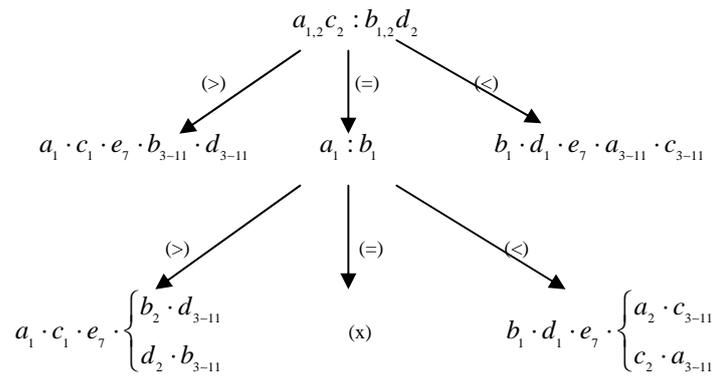

Fig.5e

The sketch of algorithm $g_1(17|_5) = 13$

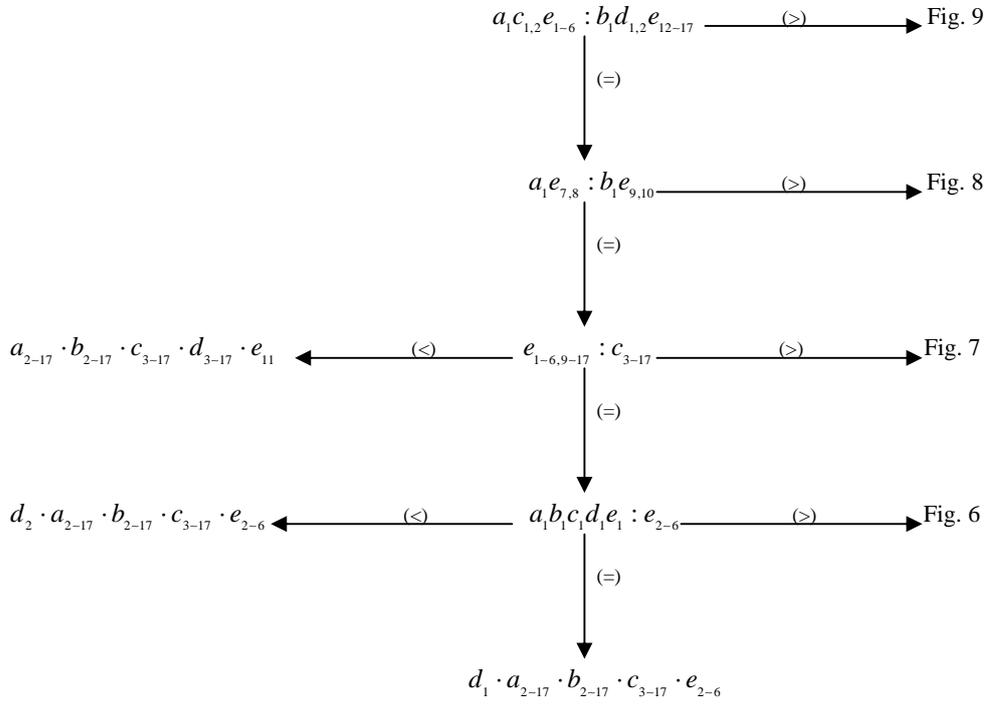

Fig. B

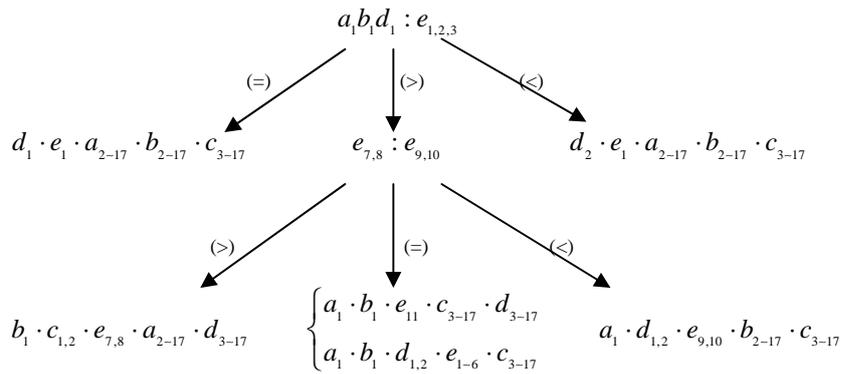

Fig. 6

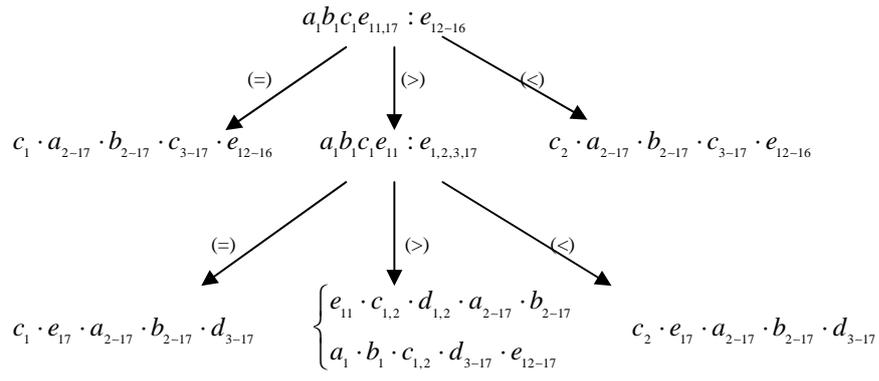

Fig. 7

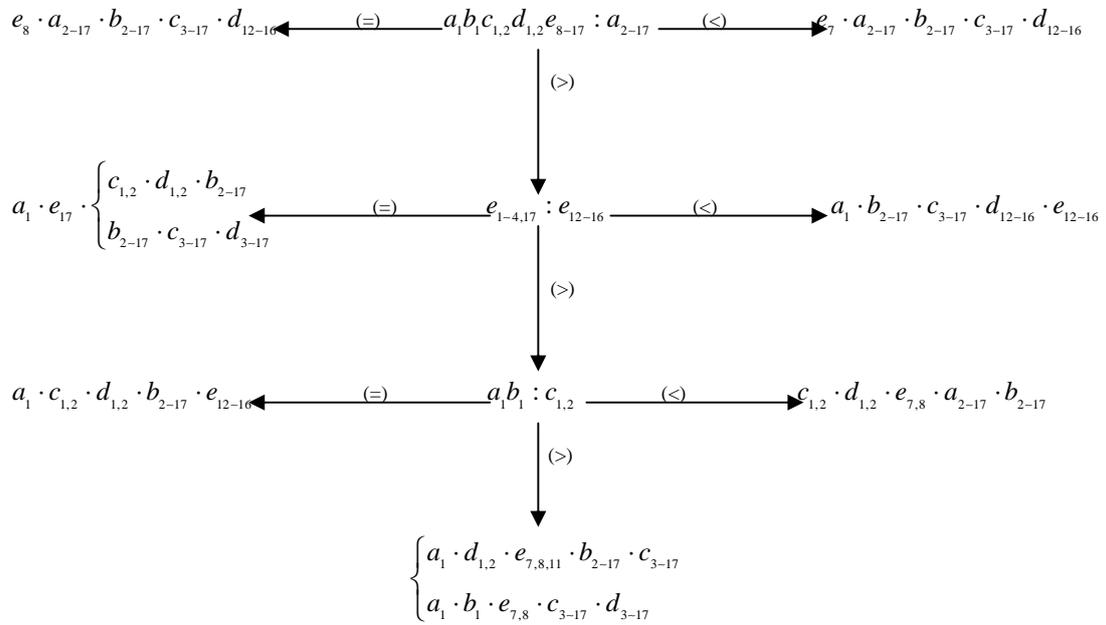

Fig. 8

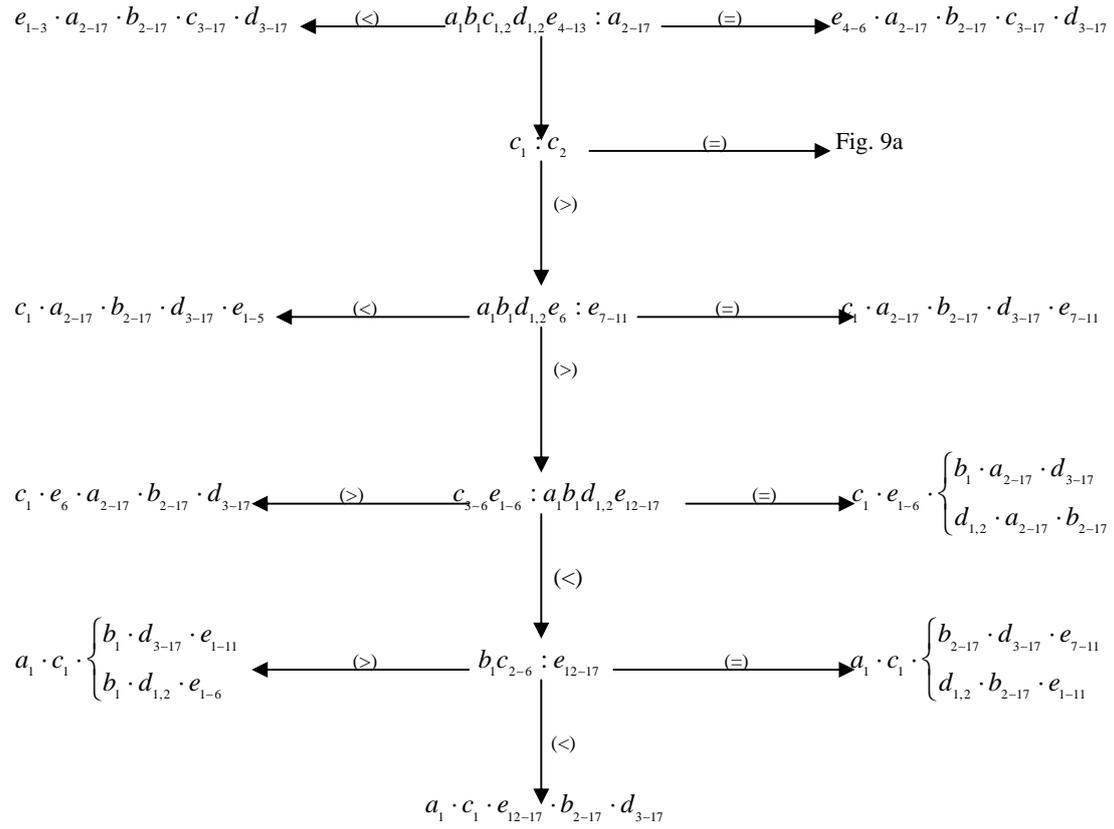

Fig. 9

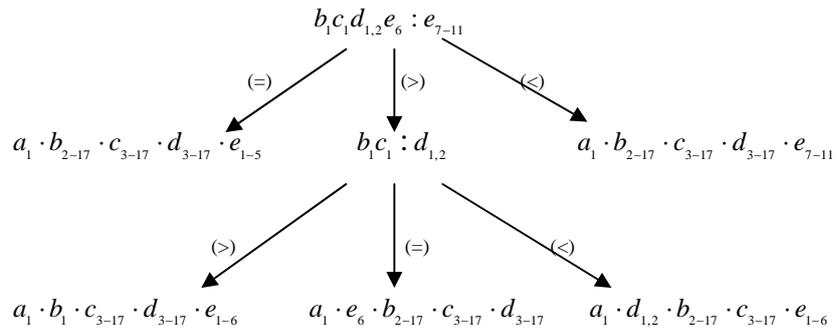

Fig. 9a

Sub-algorithm $g_1(5,16) = 4$

Let $X = \{x_i \mid 1 \le i \le 5\}$, $Y = \{y_i \mid 1 \le i \le 16\}$, $Ct(X) = Ct(Y) = 1$.

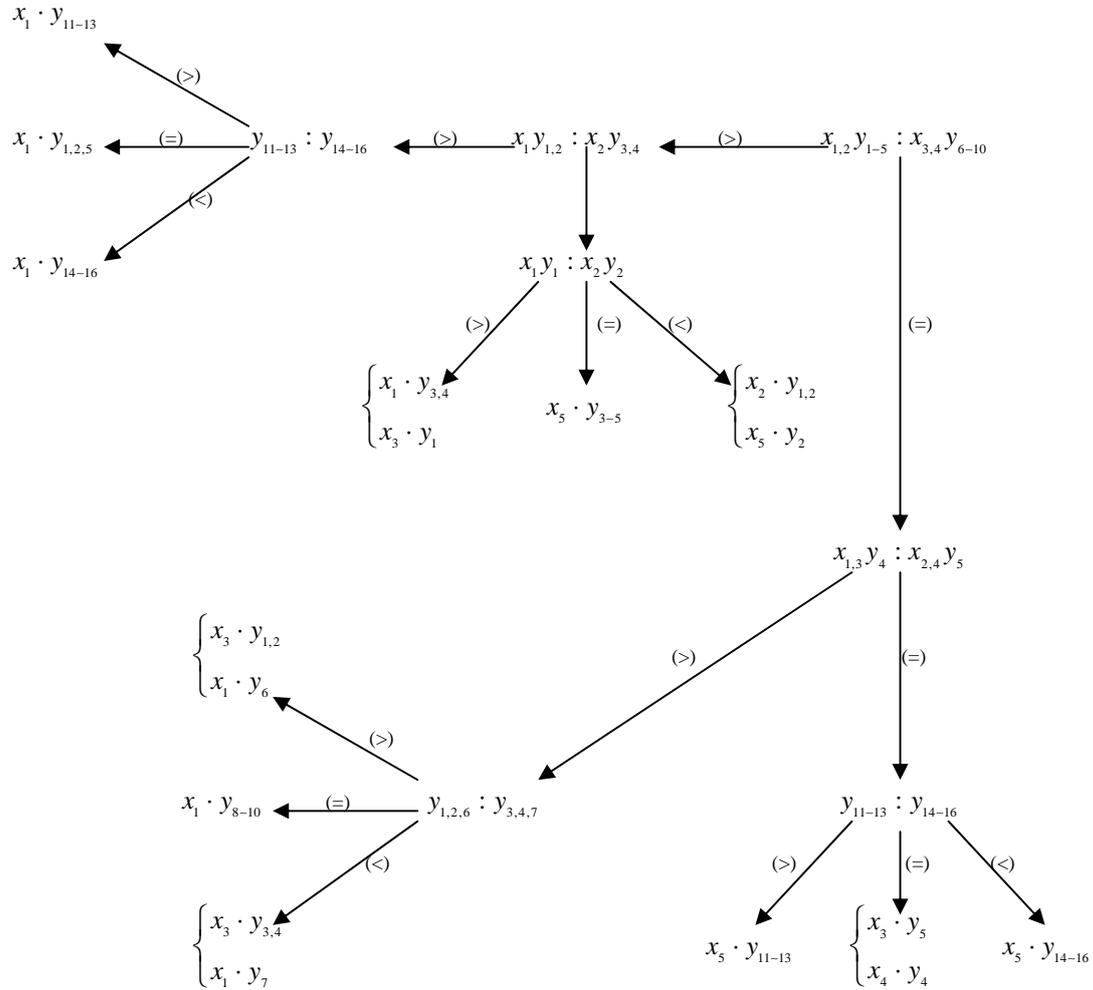

Fig. 10

The sub-algorithm $g_1(5,16) = 4$ will be used in the algorithm $g_1(17\mid_5) = 13$.

The outputs in the sketches are the objective sets settled easily by the known sub-algorithms.

The sketches of algorithms for $g_1(7\mid_5) = 9$ and $g_1(13\mid_5) = 12$ have been omitted to save the space.